\documentclass{elsart1}
\usepackage{latexsym,bm,amssymb,amsfonts,graphics,graphicx,color}

\def\RR{\mathbb R}

\def\pmatrix{ \left( \begin{array} }
\def\endpmatrix{ \end{array} \right) }
\def\aa{\alpha}

\def\bfb{\boldsymbol{b}}

\def\bfe{\boldsymbol{e}}
\def\bff{\boldsymbol{f}}

\def\bfxi{\boldsymbol{\xi}}

\def\bfpsi{\boldsymbol{\psi}}

\def\bfu{\boldsymbol{u}}
\def\bfv{\boldsymbol{v}}
\def\bfw{\boldsymbol{w}}
\def\bfx{\boldsymbol{x}}
\def\bfy{\boldsymbol{y}}

\def\bfo{\boldsymbol{0}}
\def\pmatrix{ \left( \begin{array} }
\def\endpmatrix{ \end{array} \right) }

\def\d2dxx{\frac{\partial^2}{\partial x^2}}

\def\no{\noindent}
\def\diag{{\rm diag}}
\def\proof{\underline{Proof}\quad}
\def\QED{~\mbox{$\Box$}}

\newtheorem{rmk}{Remark}

\begin{document}
\begin{frontmatter}
\title{Iterative solution of piecewise linear systems for the
numerical solution of obstacle problems\thanksref{footnoteinfo}}

\thanks[footnoteinfo]{Work developed
within the project ``Numerical methods and software for
differential equations''}

\vspace{-.5cm}
\address{Dedicated to Prof. D.\,Trigiante on the occasion of his 65$^{th}$ birthday}

\author{Luigi Brugnano}\ead{luigi.brugnano@unifi.it} ~\and
\author{Alessandra Sestini}\ead{alessandra.sestini@unifi.it}
\address{Dipartimento di Matematica ``U.\,Dini''\\
Viale Morgagni 67/A, 50134 Firenze, Italy}

\begin{abstract}
We investigate the use of {\em piecewise linear systems}, whose
coefficient matrix is a piecewise constant function of the
solution itself. Such systems arise, for example, from the numerical solution of
linear complementarity problems and in the numerical solution of
free-surface problems. In particular, we here study their
application to the numerical solution of both the (linear)
{\em parabolic obstacle problem} and the {\em obstacle problem}.
We propose a class of effective semi-iterative
Newton-type methods to find the exact solution of such piecewise
linear systems. We prove that the semi-iterative Newton-type
methods have a global monotonic convergence property, i.e., the
iterates converge monotonically to the exact solution in a finite
number of steps. Numerical examples are presented to demonstrate
the effectiveness of the proposed methods.

\end{abstract}

\end{frontmatter}

\section{Introduction}

Due to their importance in both theory and applications, since the
sixties a lot of interest has been devoted in the literature to
both the obstacle and the parabolic obstacle problems (see, e.g.,
\cite{F63,LS67}), which are nonlinear differential problems of
elliptic or parabolic type, respectively. Since the beginning,
their theoretical study has been developed within the more general
context of variational inequalities; the interested reader can
refer, e.g., to \cite{Rod}, where their mathematical-physical
introduction is given together with many related abstract
theoretical results.

In this paper we are in particular concerned with the numerical
solution of the linear obstacle and linear parabolic obstacle
problems of second order, that is we assume that the differential
operators involved in the inequalities are of second order and
that they are linear with respect to the unknown function (observe
that, despite of their name, these differential problems are still
nonlinear because they involve differential inequalities instead
of equalities). Their equivalent formulations as complementarity
problems \cite{Rod} is explicitly used in this paper for their
introduction.

The first numerical schemes for the obstacle problem were based on
finite element discretizations of such problems combined with
projected relaxation methods used for the solution of the discrete
problem \cite{Glow}. However, the convergence rate of such
iterative schemes depends on the mesh refinement, and the position
of the free boundary of the {\em coincidence set} (i.e., the set
where the solution of the problem coincides with the obstacle) is
not taken into account. In order to make the convergence rate
independent of the mesh refinement, several multigrid algorithms
have been proposed in the literature (see, e.g.,
\cite{Hackbusch,Zhang}). Another way proposed in the literature to
solve the obstacle problem is based on the use of an active set
strategy which can also be combined with the multigrid approach
(see, e.g., \cite{Hoppe1,Tarvainen}). The scheme for the solution
of the discrete problem is iterative (observe however that for the
linear case convergence is obtained in finitely many steps and it
is monotone) and at each step it makes a problem linearization by
specifying the active and the inactive part of the unknowns. The
active set strategy is used to define an outer iteration and each
inner iteration requires the solution of a reduced linear system
which can be efficiently implemented by the multigrid approach
(see, e.g., \cite{Hoppe1}) or by a preconditioned conjugate
gradient method when the linear case with Dirichelet boundary
conditions is considered  \cite{Hoppe94}. Inexact semismooth
Newton methods have been developed in \cite{Kanz}. A further
alternative is proposed in \cite{Xue}, where a different iterative
approach is considered for the linear obstacle problem. In this
case an iterative approximation of the contact region is used
(moving obstacle).

Even if in this paper we do not deal explicitly with mesh
adaptation, we consider it an important aspect to be developed in
the future in conjunction with our schemes. Relating to the coupled
active set and multilevel approach, mesh adaptation is considered
for the linear elliptic case in \cite{Hoppe94}, where some {\em a
posteriori} error estimates are reported.

Concerning the parabolic case, the use of a regularization
technique combined with a Lagrange multiplier approach is proposed
in \cite{Ito06}, where some numerical results are presented for
the Black-Scholes model for American options. In particular, such
results are obtained by using a second order (in time and space)
finite difference discretization of each regularized problem which
leads to the solution of nonsmooth nonlinear equations which are
numerically solved by a semismooth Newton method. Some interesting
numerical results related to the linear parabolic obstacle problem
are also reported in \cite{Pommier08} where an Euler implicit time
scheme is combined with a finite element spatial discretization.
In such paper some a posteriori error estimates are used in order
to control the mesh refinement, both in space and in time. Again,
the discrete problem is solved by using a semismooth Newton method
(see, e.g., \cite{Ito06}). The moving mesh method, based again on
a posteriori error estimates, is the strategy suggested in
\cite{LiMa} for both improving the accuracy and reducing the
computational cost of the finite element numerical approximations
of the parabolic obstacle problem which, otherwise, may have a
very poor efficiency. However, no numerical result is given in
that reference.

In this paper, we shall consider the numerical modeling and
solution of linear obstacle problems by means of {\em piecewise
linear systems} ({\em PLS}, hereafter), which have been recently
introduced and investigated in \cite{BC,BC1}, with application to
flows in porous media. PLS are linear systems, whose coefficient
matrix is a piecewise constant function of the solution itself.
They can be used for the efficient modeling of a number of
real-life problems. In particular, we here consider their
application for the numerical solution of the linear classical
obstacle problem and its parabolic version. The procedure proposed
for the numerical solution of the associated discrete obstacle
problems is an iteration having a monotonic finite convergence
behaviour. For completeness and clarity reasons, in the paper we
specify that such application of PLS can be also formulated as a
special case of the dual--active set strategy introduced in
\cite{Tarvainen} and there studied under the assumption that the
coefficient matrix characterizing the discrete problem is an
$M$--matrix. However, in our opinion the PLS formulation of such
iteration is an interesting alternative because of its compactness
and because it allows us to analyze the convergence features of
the method under less restrictive hypotheses which have to be
assumed when obstacle problems with Neumann boundary conditions
are dealt with.

The paper is organized as follows. In Section~\ref{clasost} we
investigate the classical obstacle problem. Then, in
Section~\ref{plss} we consider its modeling via PLS, whose
numerical solution is investigated in Sections~\ref{Newt1}. In
Section~\ref{parost} the parabolic obstacle problem is considered,
whose solution turns out to be a particular instance of what
stated in Section~\ref{Newt1}. Section~\ref{numer} contains some
numerical examples dealing with both Dirichlet and Neumann
boundary conditions and, finally, Section~\ref{fine} contains a
few conclusions.

\section{The (classical) obstacle problem}\label{clasost}
We here consider the following special linear systems which
involve nonsmooth functions of the solution itself,

\begin{equation}\label{plsnew}
\min\{\bfo,\bfx\} + T\max\{\bfo,\bfx\} = \bfb,
\end{equation}
where $\bfx = (x_i), ~\bfb = (b_i)\in\RR^n$, with $\bfb$ a known vector,

$$\max\{\bfo,\bfx\} = \pmatrix{c} \max\{0,x_1\}\\ \vdots\\
\max\{0,x_n\}\endpmatrix,\qquad \min\{\bfo,\bfx\} = \pmatrix{c}
\min\{0,x_1\}\\ \vdots\\
\min\{0,x_n\}\endpmatrix,$$ and  $T \in \RR^{n\times n}$ is a (known)
irreducible matrix, satisfying either one of the following
properties:
\begin{description}

\item[T1: ] $T$ is an $M$-matrix (i.e., it can be written as $T=\aa I -B$,
with $B\ge O$ and $\rho(B)<\aa$), or

\item[T2: ]${\tt null}(T^T) \equiv {\tt span}(\bfv)$, ${\tt null}(T)
\equiv {\tt span}(\bfw)$, ~with~ $\bfv,\bfw>\bfo$
~(componentwise), and $T+D$ is an $M$-matrix for all diagonal
matrices $D\gneqq O$ (i.e., $D\ge O$ and $D\ne O$).

\end{description}
Note that, if $\bfxi=(\xi_i)\in\RR^n$ is a given known vector,
upon a suitable variable transformation, the following problems
can be taken back to problem (\ref{plsnew}),
\begin{eqnarray}\label{eplsnew}
\min\{\bfxi,\bfx\} + T\max\{\bfxi,\bfx\} &=& \bfb,\\[2mm]
\max\{\bfxi,\bfx\} + T\min\{\bfxi,\bfx\} &=& \bfb.  \label{eplsnew1}
\end{eqnarray}

One important motivation, for solving problem (\ref{plsnew}),
stands in the efficient numerical modeling of the  the linear
obstacle problem. In more details, let us consider the problem in
its simplest form (see, e.g., \cite{Rod} for more general
formulations):
\begin{equation}\label{cop}
-\triangle u \ge f, \qquad u \ge \psi, \qquad (u-\psi)(\triangle
u+f)=0, \qquad\mbox{in}~\Omega,\end{equation}

\no with suitable prescribed boundary conditions on
$\partial\Omega$, where $f$ is a known function and $\psi$ is the
obstacle.

After a suitable finite difference discretization of problem (\ref{cop}), one
obtains a corresponding discrete complementarity problem in the form
\begin{equation}\label{copdisc} T\bfu \ge \bff, \qquad \bfu\ge \bfpsi,
\qquad (\bfu-\bfpsi)^T(T\bfu- \bff)=0, \end{equation}

\no where $\bfu$ is the unknown solution, $\bff$ depends on the
function $f$ and on the boundary conditions, $\bfpsi$ is the
discrete representation of the obstacle, and $T$ is a matrix
satisfying either {\bf T1}, if $u$ is specified in at least one
point of $\partial\Omega$, or {\bf T2},
otherwise.
the previous problem can then be reformulated as
\begin{equation}\label{lcp2}
T\bfy \ge \bfb, \qquad \bfy\ge \bfo, \qquad \bfy^T(T\bfy-
\bfb)=0,\end{equation} where $\bfb = \bff-T \bfpsi.$ The following
result then holds true.

\begin{thm}\label{thcop} If $\bfx$ is a solution of PLS
(\ref{plsnew}), then $\bfy=\max\{\bfo,\bfx\}$ is a solution of
(\ref{lcp2}).\end{thm} \proof Let $\bfx$ be a solution of
(\ref{plsnew}). Clearly, $\max\{\bfo,\bfx\}$ always satisfies the
second inequality in (\ref{lcp2}). Then, concerning the first
inequality and the complementarity condition, the following cases
can occur, when considering the generic $i$th entry of $\bfx$:
\begin{itemize}
\item $x_i<0$. Consequently, $\min\{0,x_i\}=x_i$ \,and\, $\max\{0,x_i\}=0$.
Moreover, one has that the $i$th component of  the first
inequality in (\ref{lcp2}) is satisfied. In fact, by setting
$\bfe_i$ the $i$th unit vector:
$$\bfe_i^TT\max\{\bfo,\bfx\} > \min\{0,x_i\} + \bfe_i^TT\max\{\bfo,\bfx\} =
b_i.$$

\item $x_i\ge 0$. In such a case, $\min\{0,x_i\}=0$ \,and\,
$\max\{0,x_i\}=x_i$. Moreover, the $i$th component of  the first inequality
in (\ref{lcp2}) turns out to be an equality. In fact:
$$\bfe_i^TT\max\{\bfo,\bfx\} = \min\{0,x_i\}
+\bfe_i^TT\max\{\bfo,\bfx\} = b_i.$$
\end{itemize}
Consequently, one concludes that $\bfy=\max\{\bfo,\bfx\}$
satisfies all the inequalities in (\ref{lcp2}), as well as the
complementarity condition.\QED

\section{Modeling through PLS}\label{plss}

We here show how the nonsmooth equation (\ref{plsnew}) can be
efficiently reformulated by means of a suitable PLS. In more
details, for a given vector $\bfx\in\RR^n$, let define the
following diagonal matrix:
\begin{equation}\label{Peta}
P(\bfx) =\pmatrix{ccc}p(x_1)\\ &\ddots\\ &&p(x_n)\endpmatrix,
\qquad\mbox{with}\qquad p(x_i) =\left\{ \begin{array}{cl} 1
&\mbox{~if~}
x_i\ge0,\\[2mm] 0 &\mbox{~otherwise.}\end{array}\right.\end{equation}

The following straightforward results then hold true.

\begin{lem}\label{lemma0} \quad $P(\bfx)\bfx = \max\{\bfo,\bfx\}$, \quad
$\left[I-P(\bfx)\right]\bfx = \min\{\bfo,\bfx\}$.\end{lem}

\begin{lem}\label{lemma5} System (\ref{plsnew}) is equivalent to the
following PLS:
\begin{equation}\label{plsnew1}
\left[ I-P(\bfx) + TP(\bfx)\right] \bfx = \bfb.
\end{equation}
\end{lem}

For sake of completeness, we also mention that problems
(\ref{eplsnew})--(\ref{eplsnew1}) can be respectively reformulated
as
\begin{eqnarray}
&\left[ I-P_\xi(\bfx)+ TP_\xi(\bfx)\right] (\bfx-\bfxi) = \bfb-(I+T)\bfxi,&
\label{plsnew2}\\[2mm]
&\left[ P_\xi(\bfx)+ T(I-P_\xi(\bfx))\right] (\bfx-\bfxi) = \bfb-(I+T)\bfxi,&
\label{plsnew3}
\end{eqnarray}
where
\begin{equation}\label{Pxi}
P_\xi(\bfx) =\pmatrix{ccc}\hat{p}(x_1)\\ &\ddots\\
&&\hat{p}(x_n)\endpmatrix,
\qquad\mbox{with}\qquad \hat{p}(x_i) =\left\{ \begin{array}{cl} 1
&\mbox{~if~}
x_i\ge\xi_i,\\[2mm] 0 &\mbox{~otherwise.}\end{array}\right.\end{equation}

\section{The Newton-type iteration for the obstacle
problem}\label{Newt1} Some preliminary results are stated at first
in order to derive a Newton-type procedure for solving the PLS
(\ref{plsnew1}) and prove its convergence. Their proof is
straightforward and is, therefore, omitted.

\begin{lem}\label{lemma4} Let $T$ satisfy {\bf T1}.
Then, for any diagonal matrix $P$, $O\le P\le I$, both matrices
$I-P+TP$  and $I-P+PT$ are $M$-matrices and, therefore,
$(I-P+TP)^{-1}\ge O$, \,$(I-P+PT)^{-1}\ge O$. Moreover, if in
addition $P\ne I$, the same result continues to hold when $T$
satisfies {\bf T2}.
\end{lem}

 It is to be noted that the left-hand side of system
(\ref{plsnew1}) is not everywhere differentiable.  Nevertheless, a
{\em Newton-type method} for solving system (\ref{plsnew1}) can be
deduced, $$\bfx^{k+1} = \bfx^{k} - \left( I -P^k+ TP^k
\right)^{-1}\left[ \left(
 I -P^k+ TP^k\right) \bfx^k - \bfb \right], \qquad k = 0,1,\dots,
$$
where the upper index $k$ denotes the iteration step and (see
(\ref{Peta}))
\begin{equation}\label{Pk} P^0 = O, \qquad P^k = P(\bfx^k),
\quad k=1,2,\dots.\end{equation} This simplifies to the following
Picard iteration,
\begin{equation}\label{mammanew}
 P^0 = O, \qquad \left(I-P^k+TP^k \right)\bfx^{k+1} = \bfb, \qquad
k=0,1,\dots.
 \end{equation}
The following result provides a straightforward stopping criterion
for the iteration.

\begin{lem}\label{finitonew} If, for some $k\ge 0$, one
gets \begin{equation}\label{finenew}
(P^{k+1}-P^{k})\bfx^{k+1}=\bfo,
\end{equation}
then $\bfx^*=\bfx^{k+1}$ is an exact solution of problem
(\ref{plsnew1}).
\end{lem}

\proof  Since $(P^{k+1}-P^{k})\bfx^{k+1}=\bfo$, one has
$$\left(I-P^k+TP^{k}\right)\bfx^{k+1} =
\left(I-P^{k+1}+TP^{k+1}\right)\bfx^{k+1}=\bfb,$$ i.e.,
$\bfx^{k+1}$ solves (\ref{plsnew1}).\QED\medskip

\begin{rmk}
Actually, iteration (\ref{mammanew}) combined with the stopping
criterion (\ref{finenew}) can be formulated as a special case of
the dual--active set strategy described by Algorithm $A1$ in
\cite{Tarvainen}. However, the next theorem shows that the compact
matrix formulation here considered allows us a corresponding
compact analysis of its convergence behaviour which is here
extended to the case where $T$ satisfies property ${\bf T2}$
instead of ${\bf T1}.$ In addition, we observe that the PLS
formulation (\ref{plsnew1}) of the linear discrete obstacle
problem allows us to get significant results about the existence
and uniqueness of the solution of the discrete problem even for
the extended case (see Theorem \ref{sola2}).
\end{rmk}

The iteration is well-defined under the following conditions.

\begin{thm}\label{definenew} Let matrix $T$ in system (\ref{pls1})
satisfy {\bf T1}. Then, the matrix $$\left(I-P^k+TP^k\right)$$ is
an M-matrix, and the iteration (\ref{mammanew}) is well defined
for all $k$ until convergence. If $T$ satisfies {\bf T2}, the same
result holds true, provided that
\begin{equation}\label{vb0}
\bfv^T\bfb\le0.
\end{equation}
\end{thm}

\proof The thesis easily follows from Lemma~\ref{lemma4}, if we
are able to prove that, when $T$ satisfies {\bf T2} and $P^k\ne
I$, then:
\begin{itemize}
 \item either the exit condition (\ref{finenew}) holds true,
 \item or $P^{k+1}\ne I$, so that the next iteration is well-defined.
\end{itemize}
In the first case, by virtue of Lemma~\ref{finitonew},
$\bfx^*=\bfx^{k+1}$ is a solution of the problem, so that no
further iterations are needed. In the second case, we observe
that, from the definition (\ref{Peta}), one readily shows that
$$
(P^{k+1}-P^k)\bfx^{k+1}\ge\bfo. $$ If $P^{k+1}=I$, then this would
imply that, from (\ref{vb0}) and (\ref{mammanew}), and considering
that $\bfv>\bfo$,
$$0\le \bfv^T(P^{k+1}-P^k)\bfx^{k+1} = \bfv^T(I-P^k)\bfx^{k+1} =
\bfv^T(I-P^k+TP^k)\bfx^{k+1} = \bfv^T\bfb \le0.$$ Consequently,
the exit condition (\ref{finenew}) holds true, so that
$\bfx^{k+1}$ is solution of problem (\ref{plsnew}).\QED

Next, we prove that the iteration (\ref{mammanew}) satisfies an
important property of monotony. Before that, we state the
following preliminary result, whose proof is straightforward and
is, therefore, omitted.

\begin{lem}\label{lemon} By setting as usual $\bfx^k=(x_i^k)$ and
$\bfx^{k+1}=(x_i^{k+1})$, for $k\ge1$ one has: $$\left(
P^k\bfx^{k+1}\ge P^{k-1}\bfx^k\ge \bfo \right) \Rightarrow
\left(x_i^k\ge0 ~\Rightarrow~ x_i^{k+1}\ge 0\right) \Rightarrow
\left(P^{k+1}\ge P^k\ge O\right).$$
\end{lem}

\begin{thm}\label{Pmonotona} Let the hypotheses of
Theorem~\ref{definenew} hold true. Then,
\begin{equation}\label{Pkmon}
P^{k+1}\ge P^k\ge O, \qquad k=0,1,\dots.\end{equation}
\end{thm}
\proof For $k=0$ (\ref{Pkmon}) trivially holds true, since
$P^0=O$. For $k\ge1$, let us prove, according to
Lemma~\ref{lemon}, that
\begin{equation}\label{Px}
P^k\bfx^{k+1}\ge P^{k-1}\bfx^k\ge \bfo.
\end{equation}
Since, from (\ref{mammanew}),
$$(I-P^k+TP^k)\bfx^{k+1} = (I-P^{k-1}+TP^{k-1})\bfx^k = \bfb,$$
one then obtains:
\begin{eqnarray*}
\lefteqn{\left(I-P^k+P^kT\right)P^k\bfx^{k+1} }\\
&\qquad&=~ P^k\left(I-P^k+TP^k\right)\bfx^{k+1} ~=~
P^k\left(I-P^{k-1}+TP^{k-1}\right)\bfx^k \\&&=~
\left(I-P^k+P^kT\right) P^{k-1}\bfx^k + (P^k-P^{k-1})\bfx^k.
\end{eqnarray*}
By considering that $(I-P^k+P^kT)^{-1}\ge O$,
$(P^k-P^{k-1})\bfx^k\ge \bfo$, and $P^0\bfx^1=\bfo$, (\ref{Px})
then follows.\QED

This result, allows to state the finite convergence of iteration
(\ref{mammanew}).

\begin{cor}\label{finitenew} Let T satisfy either {\bf T1} or {\bf T2}.
If $T$  satisfies {\bf T2}, assume that also (\ref{vb0}) is
satisfied. Then, iteration (\ref{mammanew}) converges in at most
$n$ steps.
\end{cor}
\proof The finite convergence easily follows from the fact that
$I\ge P^{k+1}\ge P^k\ge O$, and from the fact that, as soon as
$P^{k+1}=P^k$, then the exit condition (\ref{finenew}) is
satisfied. Obviously, by considering that $P^0=O$, this will happen
in at most $n$ steps.\QED

\begin{rem} Even though Corollary~\ref{finitenew} establishes the finite
convergence of iteration (\ref{mammanew}), nevertheless the
corresponding upper bound may be large, when the dimension of the
system is large. However, several numerical tests have shown that
convergence can occur in just a few iterates (see, e.g., the
numerical tests in Section~\ref{numer}).
\end{rem}

Next, we present a conclusion on the existence of a solution for
problem (\ref{plsnew1}). We need the following preliminary result.

\begin{lem}\label{lemma3} With reference to matrix $P$ defined in
(\ref{Peta}), for any two vectors $\bfx=(x_i)$ and $\bfy=(y_i)$,
there exists a diagonal matrix,
$$W=\diag(\omega_1, \dots, \omega_n), \qquad O\le W\le I,$$ depending on $\bfx$ and
$\bfy$, such that
\begin{equation}\label{addin2}
P(\bfx)\bfx-P(\bfy)\bfy = W \cdot \left(\bfx-\bfy\right).
\end{equation}
\end{lem}
\proof For each $i=1,2,\dots,n$, it follows from (\ref{Peta}) that
either one of the following four cases occurs:
\begin{equation}\label{casi}
\hspace{-.5cm}\begin{array}{ccccc} x_i,y_i\ge 0
&\quad\Rightarrow\quad& p(x_i)= p(y_i)=1 &\quad
\Rightarrow\quad& \omega_i=1;\\
x_i,y_i<0 &\Rightarrow& p(x_i)= p(y_i)=0  &\quad \Rightarrow\quad
&\omega_i=0;\\
x_i\ge 0> y_i &\Rightarrow& p(x_i)=1,\ p(y_i)=0 &\quad
\Rightarrow\quad &0\le \omega_i=\frac{x_i}{x_i-y_i}<1; \\
x_i< 0\le y_i &\Rightarrow& p(x_i)=0,\ p(y_i)=1 &\quad
\Rightarrow\quad& 0\le \omega_i = \frac{y_i}{y_i-x_i}<1.
\end{array}\end{equation}
This proves the validity of (\ref{addin2}). \QED\medskip

We can now state the following result.

\begin{thm}\label{sola2} Let $T$ satisfy {\bf T1}. Then a solution
to problem (\ref{plsnew1}) exists and is unique. In the case where
$T$ satisfies {\bf T2}, then a solution of problem
(\ref{plsnew1}):
\begin{itemize}
\item exists and is unique when $\bfv^T\bfb<0$;

\item exists but is not unique when $\bfv^T\bfb=0$;

\item doesn't exist when $\bfv^T\bfb>0$.

\end{itemize}
\end{thm}

\proof Concerning the existence of a solution, the thesis follows
from Corollary~\ref{finitenew}, when $T$ satisfies {\bf T1}, or
$T$ satisfies {\bf T2} and (\ref{vb0}) holds true. It remains to
prove that no solution exists when $T$ satisfies {\bf T2} and
$\bfv^T\bfb>0$. Indeed, if such a vector $\bfx$ would exist, by
considering that $\bfv>\bfo$ and taking into account
(\ref{plsnew1}), then $$0<\bfv^T\bfb =
\bfv^T[I-P(\bfx)+TP(\bfx)]\bfx = \bfv^T[I-P(\bfx)]\bfx\le0,$$
which is clearly impossible.

Concerning uniqueness, let $\bfx$ and $\bfy$ be two solutions of
(\ref{plsnew1}). Then $$\bfb = [I-P(\bfx)+TP(\bfx)]\bfx =
[I-P(\bfy)+TP(\bfy)]\bfy.$$ By virtue of Lemma~\ref{lemma3}, this
implies that $$M\cdot(\bfx-\bfy) \equiv (I-W+TW)(\bfx-\bfy) =
\bfo,$$ for a suitable diagonal matrix $W$, $O\le W\le I$.
Consequently, by taking into account Lemma~\ref{lemma4}:
\begin{itemize}
\item if $T$ satisfies {\bf T1}, then $M$ is nonsingular
 and uniqueness $(\bfx=\bfy)$ follows;

\smallskip
\item if $T$ satisfies {\bf T2}, then $M$ is nonsingular,
 and uniqueness $(\bfx=\bfy)$ follows, if and only if $W\ne I$. By
 considering the possible cases (\ref{casi}), this is
 equivalent to requiring that at least one entry of one the two vectors $\bfx$
 and $\bfy$ is negative. This is indeed the case, when
 $\bfv^T\bfb<0$, since
 $$0>\bfv^T\bfb = \bfv^T[I-P(\bfx)+TP(\bfx)]\bfx = \bfv^T[I-P(\bfx)]\bfx,$$ which, by
 considering that $\bfv>\bfo$ and (see
Lemma~\ref{lemma0}) $[I-P(\bfx)]\bfx\le\bfo$,
 implies that at least one entry of $\bfx$ is negative. On the
 other hand, when $\bfv^T\bfb=0$, then $$0=\bfv^T\bfb =
 \bfv^T[I-P(\bfx)]\bfx=\bfv^T[I-P(\bfy)]\bfy,$$ which implies that
 $P(\bfx)=P(\bfy)=I$. Consequently, one obtains $T\bfx=T\bfy$, i.e., $$\bfx-\bfy\in{\rm null}(T)
 \equiv{\rm span}(\bfw).$$ In  more details,  if $\bfx=(x_i)$ is a solution
 of (\ref{plsnew1}) such that $$\min_i x_i = 0,$$ (observe that, since $\bfw>\bfo$,
 such a solution always exists), then all solutions of (\ref{plsnew1})
 are given by $$\bfx(\aa) = \bfx + \aa\bfw, \qquad \aa\ge0.\QED$$
\end{itemize}

\begin{rem} It is remarkable to observe that iteration
(\ref{mammanew}) converges to a solution of problem
(\ref{plsnew1}) under the same hypotheses that guarantee its
existence. Moreover, convergence to a solution is guaranteed also
when there is no uniqueness (i.e., when matrix $T$ satisfies {\bf
T2} and $\bfv^T\bfb=0$).
\end{rem}

For completeness, we mention that for solving problem
(\ref{eplsnew}), i.e. for solving the PLS (\ref{plsnew2}), the
corresponding iteration is:
$$P_\xi^0 = O, \qquad \left(I-P_\xi^k+T P_\xi^k\right)
(\bfx^{k+1}-\bfxi)=\bfb-(I+T)\bfxi, \qquad k=0,1,\dots,$$ where,
see (\ref{Pxi}), $P_\xi^k = P_\xi(\bfx^k)$. The same iteration can
be used for solving problem (\ref{eplsnew1}), i.e. for solving the
PLS (\ref{plsnew3}), if $P_{\xi}^k$ is replaced with $\hat
P_{\xi}^k = (I-P_{\xi}^k).$

\section{The parabolic obstacle problem}\label{parost}
We now consider the application of PLS for the numerical solution
of special linear systems, involving nonsmooth functions of the
solution itself, in the form

\begin{equation}\label{pls}
\bfx + T\max\{\bfo,\bfx\} = \bfb,
\end{equation}
where, as before, matrix $T$ satisfies either {\bf T1} or {\bf
T2}. One important motivation, for solving problem (\ref{pls}),
stands in the efficient numerical modeling of the linear parabolic
obstacle problem. In more details, let us consider the problem in
its simplest form (see, e.g., \cite{PeSha} for more general
formulations):
\begin{eqnarray}\nonumber
&&u_t \ge \triangle u +f, \qquad u \ge \psi, \\
\label{pop}\\ \nonumber &&(u-\psi)(u_t-\triangle u-f)=0,
\qquad\mbox{in}~\Omega, \qquad\mbox{for}~t>0,\end{eqnarray}

\no with suitable prescribed initial and boundary conditions at
$t=0$ and on $\partial\Omega$. Here $u_t$ is the partial time
derivative of the (unknown) solution $u$, $f$ is a known function,
and $\psi$ is the (known) function describing the {\em obstacle}.
A suitable implicit, finite difference discretization of problem
(\ref{pop}), generates a corresponding discrete complementarity
problem in the form
$$\bfu^{n+1} +T\bfu^{n+1} \ge \bfu^n +\bff,  \qquad \bfu^{n+1}\ge \bfpsi,$$
$$(\bfu^{n+1}-\bfpsi)^T(\bfu^{n+1} +T\bfu^{n+1}- \bfu^n-\bff)=0,$$

\no where $\bfu^n$ is the discrete approximation at the $n$th time
step ($\bfu^0$ being specified by the initial condition), the
vector $\bff$ depends on the function $f,$ on the boundary
conditions and on the timestep, $\bfpsi$ is the discrete
representation of the obstacle, and $T$ is a matrix satisfying
either {\bf T1} if, in the boundary conditions at the $(n+1)$st
time step, $u$ is specified in at least one point of
$\partial\Omega$, or {\bf T2}, otherwise.\footnote{As in the
previous case, the problems in the literature generally prescribe
the value of the solution at the boundary. For completeness,
however, also in this case we consider the more general kind of
boundary conditions.} By setting $\bfy=\bfu^{n+1}-\bfpsi$ and by
defining a suitable (known) vector $\bfb$, the previous problem,
to be solved at each time step, can be reformulated as
\begin{equation}\label{lcp1}
\bfy +T\bfy \ge \bfb, \qquad \bfy\ge \bfo, \qquad \bfy^T(\bfy
+T\bfy- \bfb)=0.\end{equation}

The following result then holds true.
\begin{thm}\label{thpop} If $\bfx$ is a solution of (\ref{pls}), then
$\bfy=\max\{\bfo,\bfx\}$ is a solution of (\ref{lcp1}).\end{thm}

\proof Let $\bfx$ be a solution of (\ref{pls}). Clearly,
$\max\{\bfo,\bfx\}$ always satisfies the second inequality in
(\ref{lcp1}). Then, concerning the first inequality and the
complementarity condition, the following cases can occur, when
considering the generic $i$th entry of $\bfx$:
\begin{itemize}
\item $x_i<0$. Consequently, $\max\{0,x_i\}=0$. Moreover, one has that the
$i$th component of the first inequality in (\ref{lcp1}) is
satisfied. Indeed, by setting $\bfe_i$ the $i$th unit vector:$$
\max\{0,x_i\} + \bfe_i^TT\max\{\bfo,\bfx\} > x_i +
\bfe_i^TT\max\{\bfo,\bfx\} = b_i.$$

\item $x_i\ge 0$. In such a case, $\max\{0,x_i\}=x_i$. Moreover, the $i$th
component of the first inequality in (\ref{lcp1})  turns out to be
an equality. In fact: $$\max\{0,x_i\} + \bfe_i^TT\max\{\bfo,\bfx\}
= x_i + \bfe_i^TT\max\{\bfo,\bfx\} = b_i.$$
\end{itemize}
One then concludes that $\bfy=\max\{\bfo,\bfx\}$ satisfies all the
inequalities in (\ref{lcp1}), as well as the complementarity
condition.\QED

From Lemma~\ref{lemma0}, the following straightforward result
follows.

\begin{lem}\label{lemma2} System (\ref{pls}) is equivalent to the
following PLS:
\begin{equation}\label{pls1}
\left[ I+ TP(\bfx)\right] \bfx = \bfb.
\end{equation}
\end{lem}

We now present an iterative procedure for solving (\ref{pls1}),
whose analysis will be taken back to what stated in
Section~\ref{Newt1}. By using similar arguments to those used in
that section, the iteration for solving (\ref{pls1}) is given by:
\begin{equation}\label{mamma}
 P^0 = O, \qquad \left(I+TP^k \right)\bfx^{k+1} = \bfb, \qquad
k=0,1,\dots,\end{equation} where $P^k$ is given, as usual, by
(\ref{Pk}) and matrix $T$ satisfies either {\bf T1} of {\bf T2}.
However, we observe that iteration (\ref{mamma}) can be formally
rewritten as
$$
 P^0 = O, \qquad \left[I-P^k+(I+T)P^k \right]\bfx^{k+1} = \bfb, \qquad
k=0,1,\dots.$$ This implies that, since matrix $(I+T)$ is
obviously {\bf T1}, for the iteration (\ref{mamma}) hold all the
results stated in the previous Section~\ref{Newt1} for the
iteration (\ref{mammanew}), in the {\bf T1} case. Namely, the
iteration (\ref{mamma}) always converges to the (unique) solution
of problem (\ref{pls1}) in a finite number of steps.

\section{Numerical tests}\label{numer}

In this section, we report a few numerical results for the above
iterative methods. In particular, in Section~\ref{copnum} we test
the iteration (\ref{mammanew}) for solving the classical obstacle
problem, whereas in Section~\ref{popnum} we test the iteration
(\ref{mamma}) for the parabolic obstacle problem. As it seems more
frequent in the literature, we mainly consider problems with
Dirichlet boundary conditions and, as a consequence, discrete
problems with matrices $T$ satisfying property {\bf T1}. However,
in order to present also some results related to the extended case
of a matrix $T$ satisfying {\bf T2}, a variant with Neumann
boundary conditions is also examined, for both problems presented
in the elliptic case.

We specify  that in all the reported tests, the linear systems
associated with our iteration, which are sparse and nonsymmetric,
are solved by using QMR (which is a standard iterative solver for
unsymmetric systems; see, e.g., \cite{Saad}).\footnote{We mention
this only for sake of completeness: actually, it is not intended
to discuss here the efficient solution of such linear systems
which, by the way, can be further improved by using a suitable
preconditioning technique.}

\subsection{The obstacle problem}\label{copnum}
\begin{figure}[p]
\centerline{\includegraphics[height=.4\textheight]{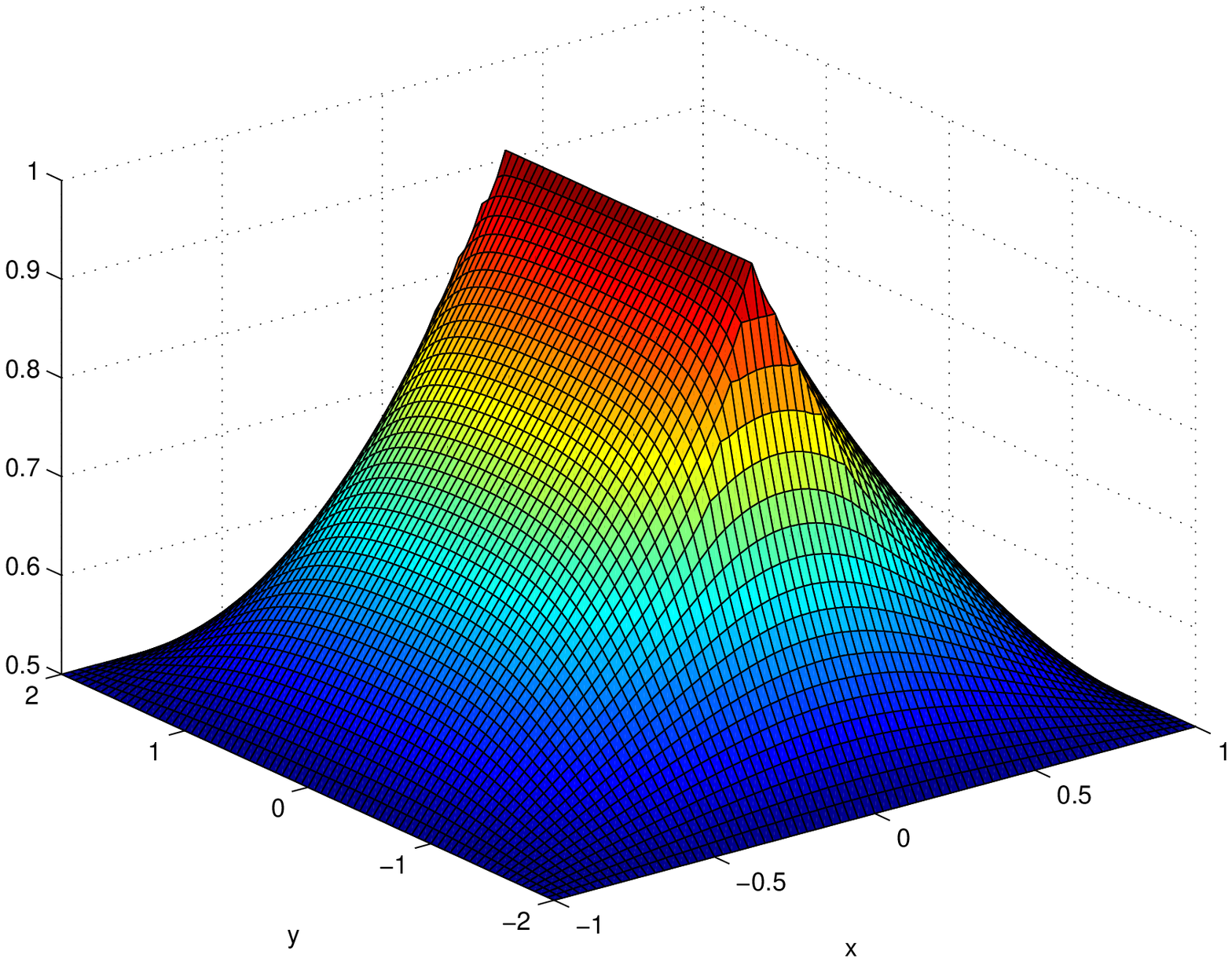} }
\caption{\protect\label{tentf} Solution of problem (\ref{tent}).}
\bigskip
\centerline{\includegraphics[height=.4\textheight]{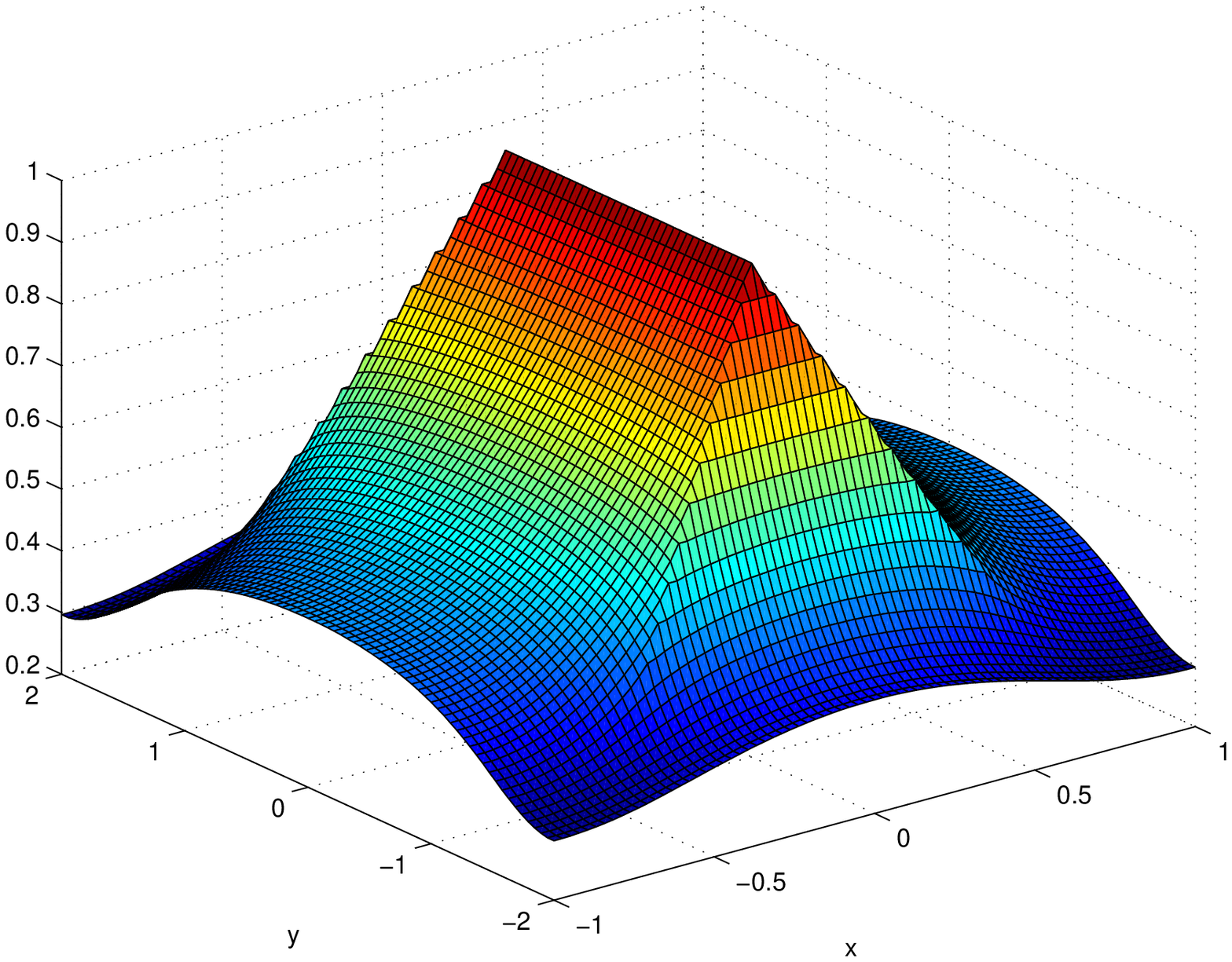}
} \caption{\protect\label{tentf-var} Solution of the variant of
problem (\ref{tent}) with homogeneous Neumann boundary
conditions.}
\end{figure}

Let us consider the following problem, whose obstacle is a non-smooth
``tent-shaped'' function:
\begin{eqnarray}\nonumber
&&-\triangle u \ge 0, \qquad u(x,y)\ge \min(1-|x|,2-|y|)\equiv \psi(x,y), \\
\label{tent}\\ && \Delta u(u-\psi)=0, \qquad (x,y)\in\Omega =
(-1,1)\times(-2,2), \qquad
u\mid_{\partial\Omega} \equiv \frac{1}2.\nonumber\end{eqnarray}
Problem (\ref{tent}) is discretized, by using the standard
5-points second-order difference scheme, on a cartesian grid with
stepsizes
\begin{equation}\label{passi}
2\Delta x = \Delta y = \frac{4}{N+1}.
\end{equation}
The resulting PLS has then dimension $n = N^2$. Because of the
Dirichlet boundary conditions, the matrix $T$ of such PLS, see
Section~\ref{clasost}, turns out to satisfy {\bf T1}.
Figure~\ref{tentf} shows the plot of the computed numerical
solution, whereas in Table~\ref{tentt} we list the number of
iterations of (\ref{mammanew}), $K$, required to get convergence,
for various values of $N$.

A variant of the above problem is also considered. In this case we
assume homogeneous Neumann boundary conditions and a non vanishing
forcing function $f$ constantly equal to $-1$.\footnote{The last
change has been introduced in order to deal with a discrete
problem admitting a unique solution.} The Neumann boundary
conditions have been discretized by using the standard 3-points
second order forward an backward difference scheme. After
eliminating the boundary unknowns, the associated PLS has always
dimension $n=N^2$ but in this case the corresponding matrix $T$
satisfies {\bf T2} and ${\bf v} = {\bf w}$ with all unit
components. The forcing term $f$ is such that the right--hand side
vector $\bfb$ in (\ref{plsnew1}) has negative sum and, thus,
Theorem \ref{sola2} allows us to state that the discrete problem
admits a unique solution. Figure~\ref{tentf-var} shows the plot of
the computed numerical solution, whereas in Table~\ref{tentt} we
list the number of iterations of (\ref{mammanew}), $K_V$, required
to get convergence, for various values of $N$.

\begin{figure}
\bigskip
\centerline{\includegraphics[height=.4\textheight]{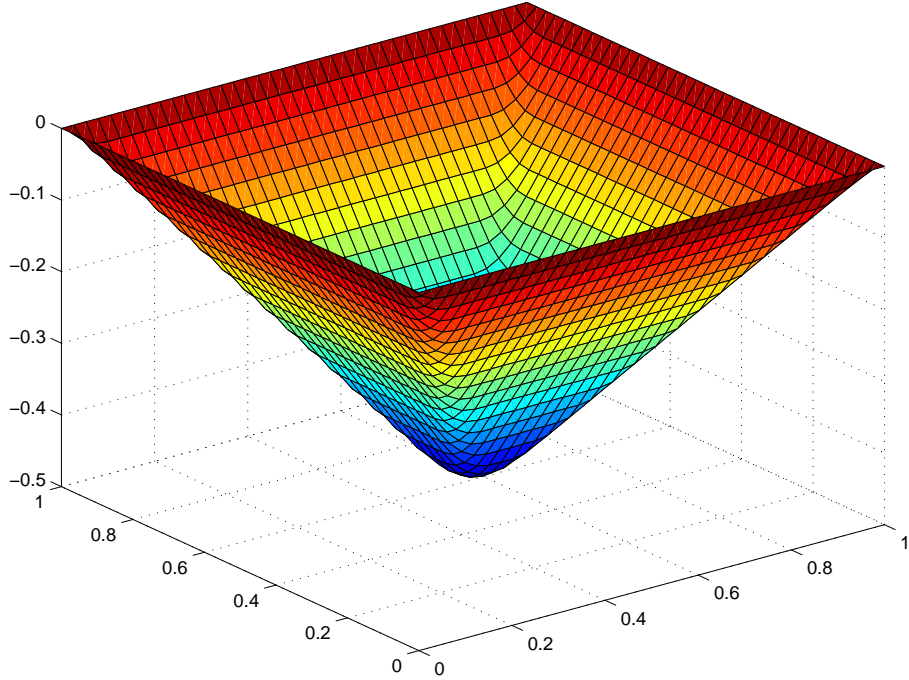}}
\caption{\protect\label{eppf} Solution of problem (\ref{epp}) with
Dirichlet boundary conditions (\ref{Dir-elasto}), $C=-20$.}
\end{figure}

\begin{figure}
\centerline{\includegraphics[width=6.5cm,height=6.5cm]{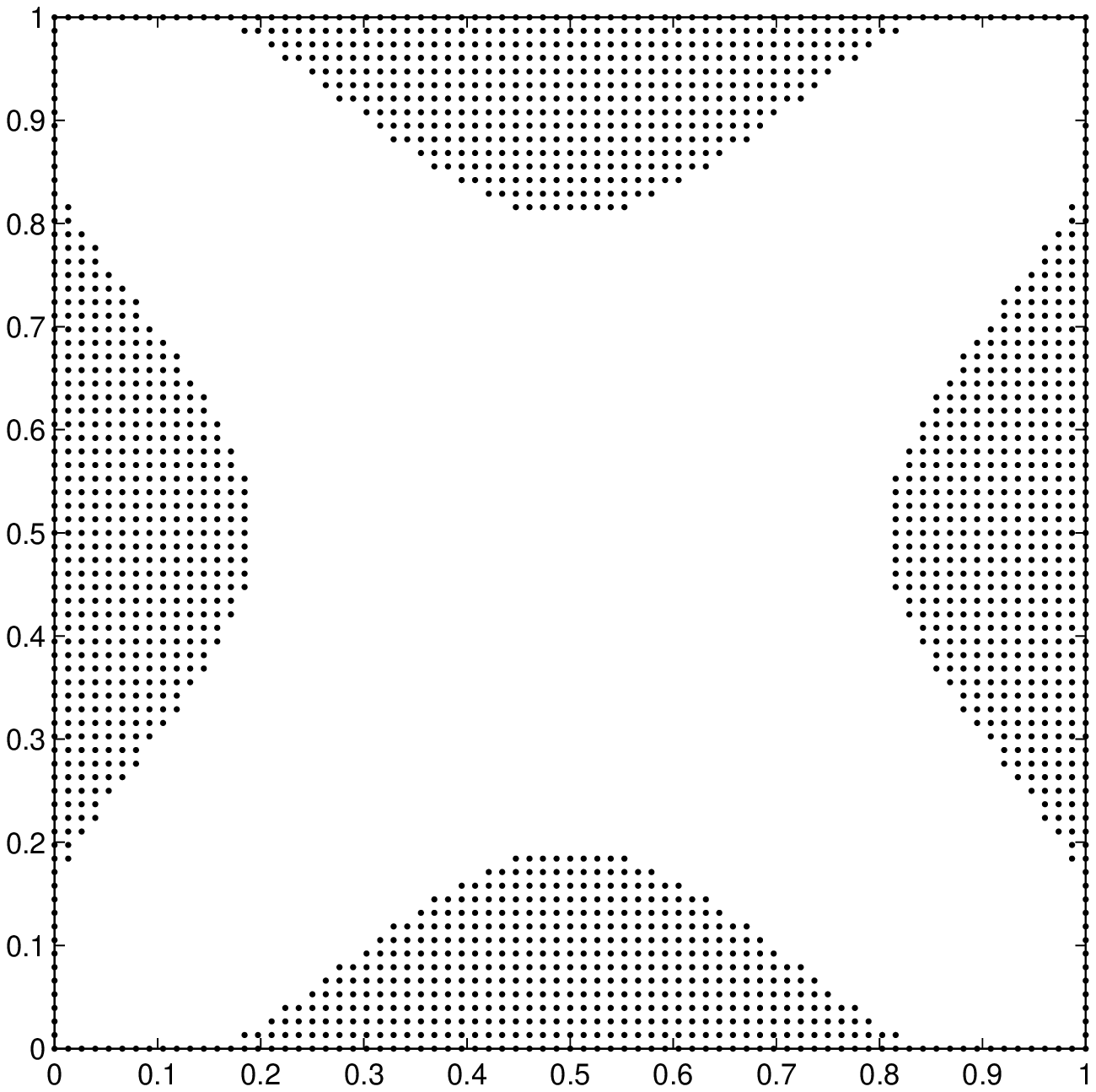}
\includegraphics[width=6.5cm,height=6.5cm]{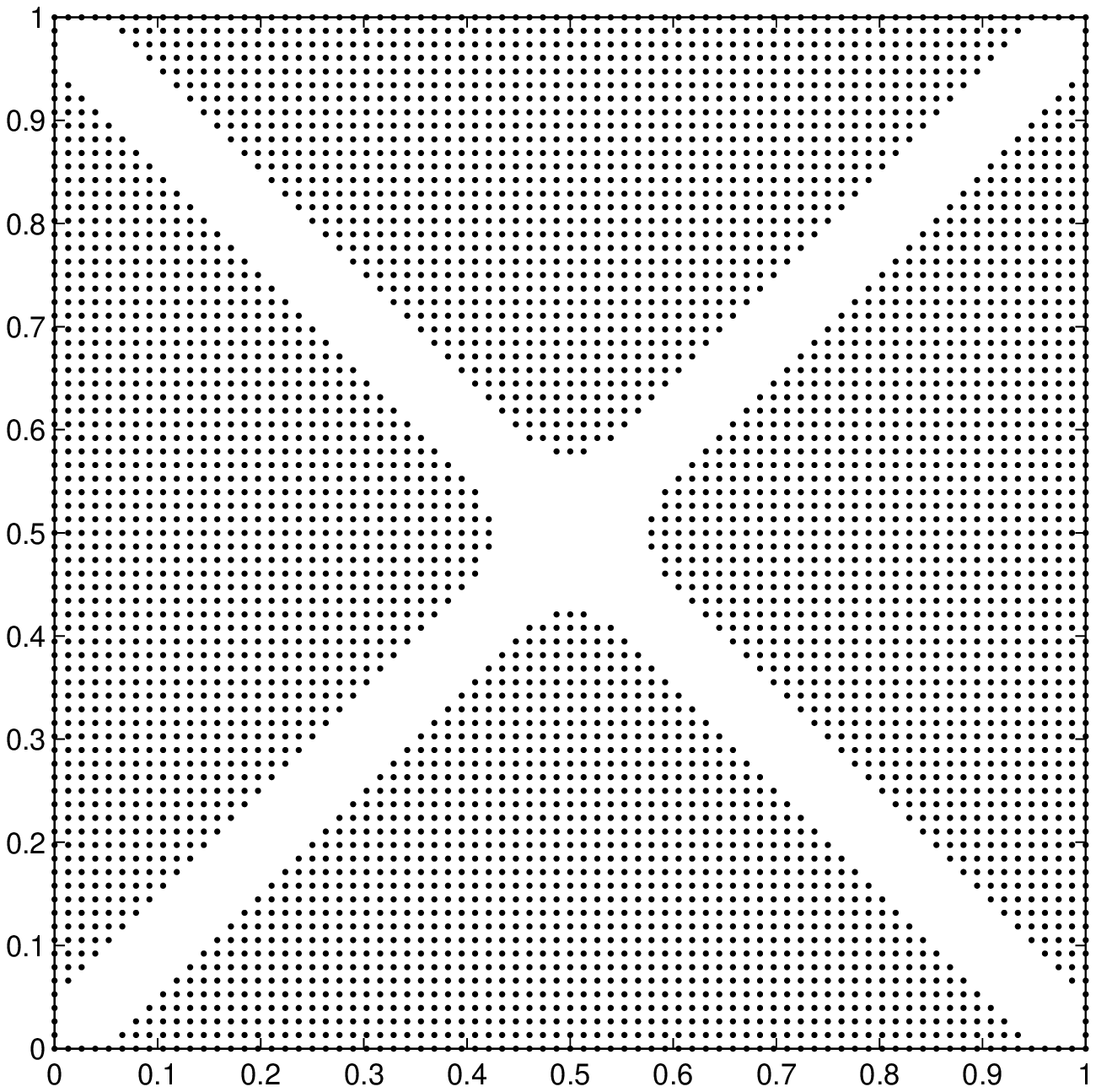}}
\caption{\protect\label{C_elastoD} Coincidence set (dotted region)
related to problem (\ref{epp}) with Dirichlet boundary conditions
(\ref{Dir-elasto}). On the left for $C=-5$ and on the right for
$C=-20.$}
\centerline{\includegraphics[width=6.5cm,height=6.5cm]{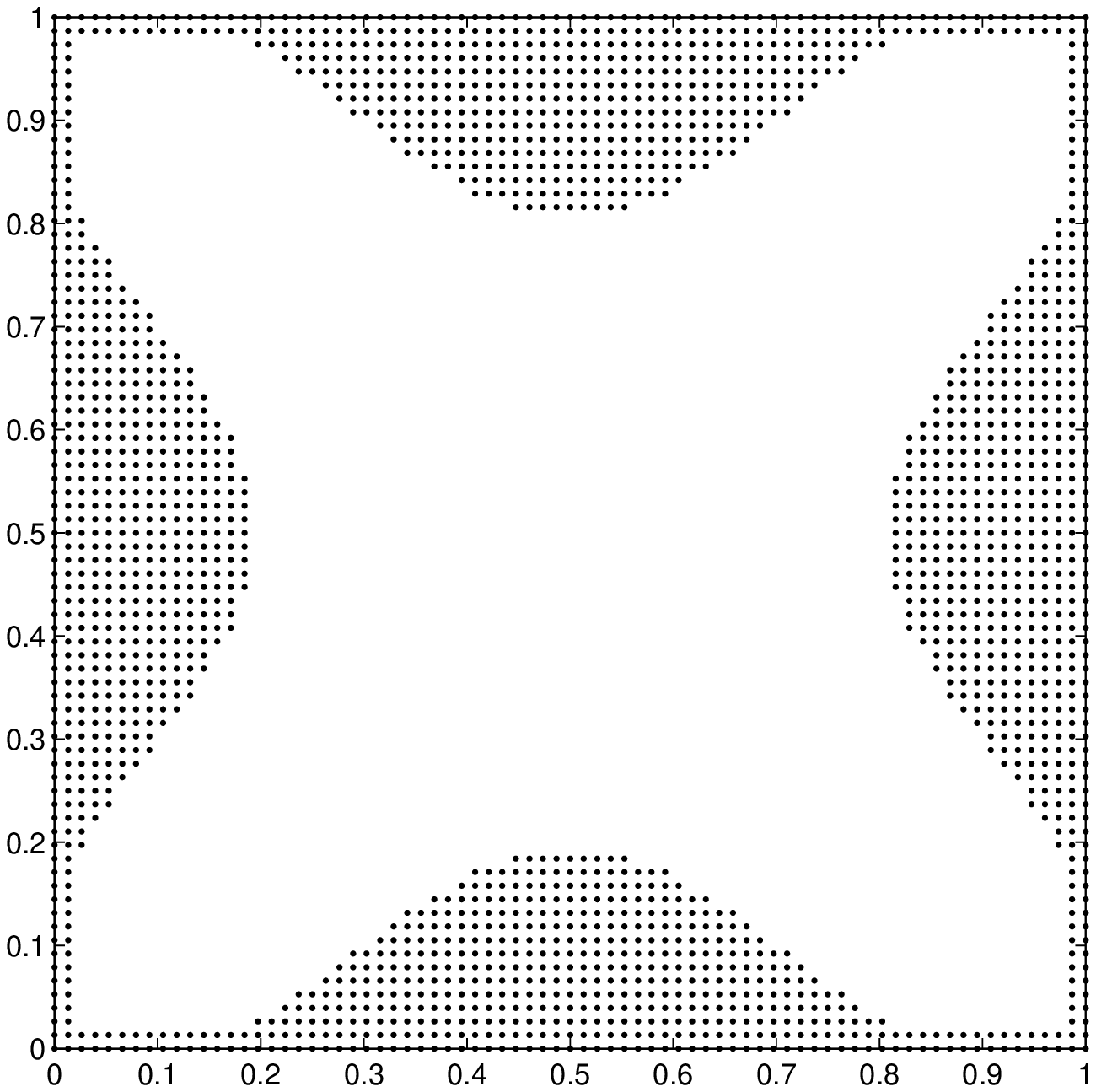}
\includegraphics[width=6.5cm,height=6.5cm]{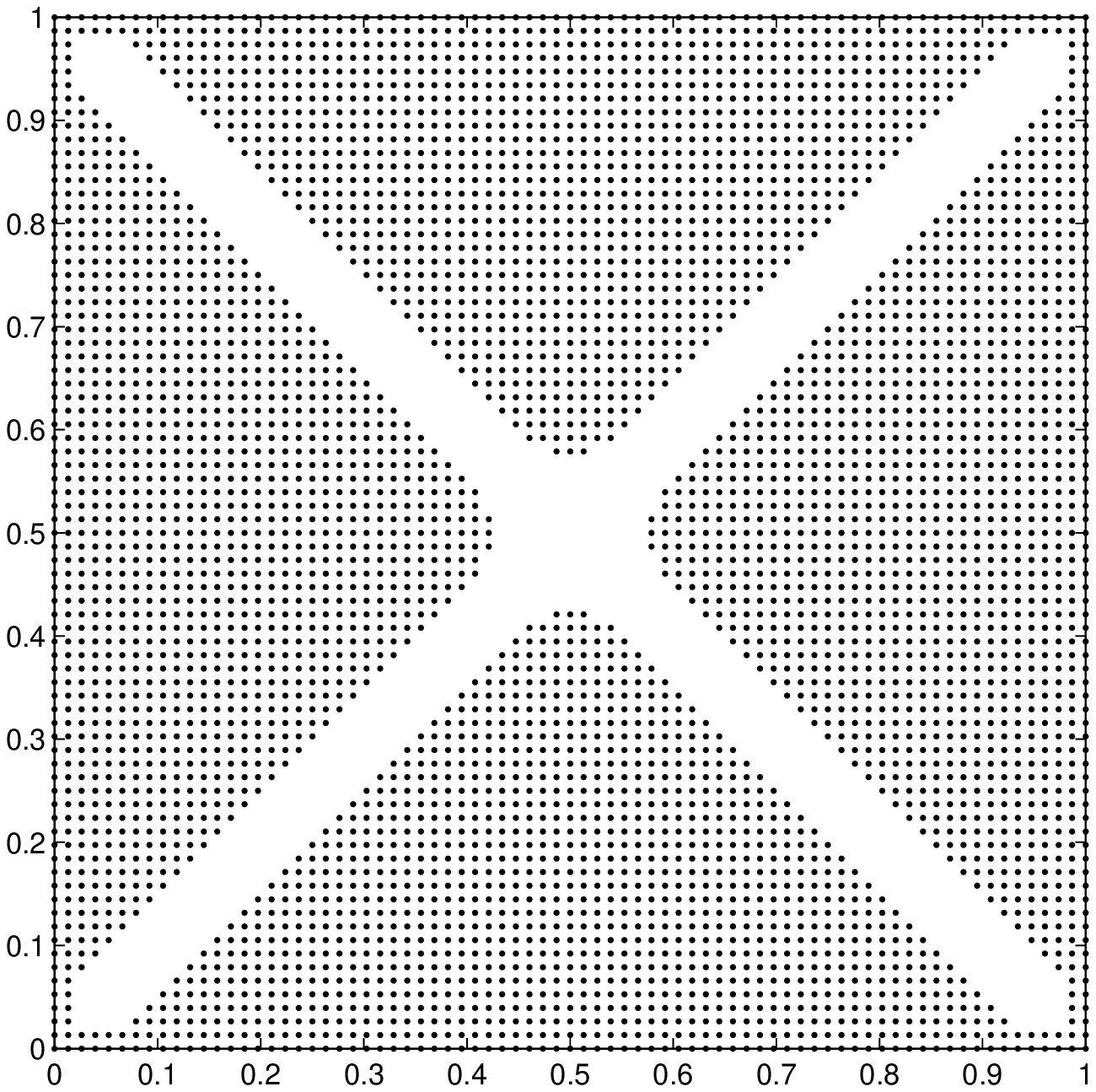}}
\caption{\protect\label{C_elastoN} Coincidence set (dotted region)
related to problem (\ref{epp}) with Neumann boundary conditions
(\ref{Neum-elasto}). On the left for $C=-5$ and on the right for
$C=-20.$}
\bigskip
\end{figure}

The second test problem is the elastic-plastic torsion problem in
\cite{Xue},
\begin{equation}
\label{epp}
\begin{array}{ll}
&-\triangle u \ge C, \qquad u(x,y)\ge
-\min(x,1-x,y,1-y)\equiv\psi(x,y), \cr
 &(\Delta u +C)(u-\psi)=0,
\qquad (x,y)\in\Omega = (0,1)^2,  \cr
 \end{array}
\end{equation}
  where $C<0$ is a given
constant, and with homogeneous \break Dirichlet boundary
conditions,
\begin{equation}\label{Dir-elasto} u \mid_{\partial\Omega} \equiv
\psi\mid_{\partial\Omega} \equiv 0.
\end{equation}
It is known that the larger $|C|$, the more difficult the problem.
Also in this case, a standard discretization, on a cartesian grid
with stepsizes
\begin{equation}\label{passi1}
\Delta x = \Delta y = \frac{1}{N+1},
\end{equation}
leads to a PLS of dimension $n=N^2$, whose matrix $T$ satisfies
{\bf T1}. Figure~\ref{eppf} shows the plot of the computed
numerical solution ($C=-20$), whereas in Table~\ref{eppt} we list
the number of iterations, $K(C)$, required to get convergence for
different values of $C$ and of the discretization parameter $N$ in
(\ref{passi1}). In such a case, the number of the required
iterations of (\ref{mamma}) decreases, as $|C|$ increases. This
behaviour can be explained considering that, as $|C|$ increases,
the coincidence set, i.e. the set of points in $\Omega$ where $u =
\psi,$ enlarges, as shown in Figure \ref{C_elastoD}. Consequently,
our initialization becomes nearer to the solution. In fact, in our
iterative procedure (\ref{mammanew}) we assume $P^0 = O;$ thus
$\bfx^1 = \bfb $ and, for sufficiently fine grids, the negativity
of $C$ and the analytical expression of the obstacle function
imply that all, or almost all, the components of $\bfb $ are
negative. Thus, as the solution $\bfu$ of (\ref{copdisc}) is
initialized with $\max\{ {\bf 0},\bfx^1 \} +\bfpsi,$ its
initialization in this case is almost everywhere coincident with
the obstacle and, thus, closer and closer to the final solution as
$|C|$ is increased.

\begin{table}[t]
\caption{\protect\label{tentt} Numerical results for problem
(\ref{tent}) and its Neumann variant, both discretized with
stepsizes (\ref{passi}); $K$ and $K_V$ respectively denote the
number of iterations of (\ref{mammanew}) for problem (\ref{tent})
and for its variant.} \centerline{\begin{tabular}{||c||r|r|r|r||}
\hline\hline $N$ & 25 & 50 & 75 & 100\\ \hline \hline $n$ &  625 &
2500 & 5625 & 10000
\\ \hline \hline
$K$ &  6 & 10 & 10 & 12\\
 \hline \hline
$K_V$ &  12 & 25 & 37 & 49\\
\hline\hline
\end{tabular}}
\bigskip
\caption{\protect\label{eppt} Numerical results for problem
(\ref{epp}) discretized with stepsizes (\ref{passi1}), either with
Dirichlet boundary conditions (\ref{Dir-elasto}) or with Neumann
boundary conditions (\ref{Neum-elasto}); $K(C)$ denotes the number
of iterations of (\ref{mammanew}) for the specified value of the
parameter $C$.} \centerline{\begin{tabular}{||c||r|r|r|r||}
\hline\hline $N$ & 25 & 50 & 75 & 100\\ \hline \hline $n$ &  625 &
2500 & 5625 & 10000
\\ \hline \hline
$K(C=-5)$  &  9 & 17 & 25 & 32\\
$K(C=-10)$ &  5 & 10 & 13 & 16\\
$K(C=-15)$ &  4 & 7 & 9 & 11\\
$K(C=-20)$ &  4 & 5 & 7 & 9\\
\hline\hline
\end{tabular}}
\end{table}
Even for this problem we have considered a variant with Neumann
boundary conditions, which have been chosen in order to get a
solution with a shape analogous to that of the solution of the
Dirichlet problem. In more detail, (\ref{Dir-elasto}) is replaced
by the following non homogeneous Neumann boundary conditions,
\begin{equation} \label{Neum-elasto}
\left. \frac{\partial u}{\partial n} \right| _{
\partial\Omega} = \left. \frac{\partial \psi}{\partial n}
\right|_{\partial\Omega},
\end{equation}
where a suitable extension of the derivative of $\psi$ at the
corner points of the domain is used. By using again a second order
discretization of the boundary conditions and eliminating the
boundary unknowns, the problem dimension is unchanged and, as for
the variant of the first problem, we deal with a matrix $T$
satisfying {\bf T2}, and ${\bf v} = {\bf w}$ with all unit
components. Also in this case, from Theorem \ref{sola2}, we obtain
a unique solution for the associated discrete problem. For all the
values of $N$ and $C$ considered in Table \ref{eppt}, as outlined
in the caption of the table, convergence has been obtained with
the same number of iterations as for the Dirichlet case. No plot
of the solution is reported for this problem with Neumann
conditions because, for all the considered values of $C,$ it is
analogous to that related to the Dirichlet case. We only present
the obtained coincidence sets in Figure \ref{C_elastoN} which,
when compared with the corresponding ones in Figure
\ref{C_elastoD}, better show the small difference near the
boundary between the obtained solution and that related to the
Dirichlet case.

\subsection{The parabolic obstacle problem}\label{popnum}

The problems that we shall consider here, are evolutionary versions of those
considered in Section~\ref{copnum}. In more details, the first problem is given
by
\begin{eqnarray}\nonumber
&&u_t \ge \triangle u, \qquad u(x,y,t)\ge \psi(x,y), \qquad (u_t-\Delta
u)(u-\psi)=0,\\
\label{ptent}\\ && (x,y,t)\in\Omega\times(0,\tau], \qquad
u\mid_{\partial\Omega} \equiv \frac{1}2, \qquad u\mid_{t=0} =
\max\left(\psi,\frac{1}2\right),\nonumber\end{eqnarray} where
$\psi$ and $\Omega$ are the same items defined in (\ref{tent}).
The spatial discretization is the same used for that problem (see
(\ref{passi})), whereas the discretization in time is, for sake of
simplicity, done by means of the implicit Euler method, by using a
constant stepsize
\begin{equation}\label{Dt} \Delta t = \frac{\tau}{\nu},\end{equation}
being $\nu$ the number of time steps. The resulting PLS, to be
solved at each time step, has dimension $n=N^2$, whose matrix $T$
is the same as that obtained for problem (\ref{tent}).
Table~\ref{ptentt} summarizes the obtained results, in terms of
required number of iterations, for $\tau=10^4$ and $\nu=20$. In
such a case, the approximation at $t=\tau$ is quite close to the
limit solution plotted in Figure~\ref{tentf}. Here the number of
iterations of (\ref{mamma}) appears to be independednt of spatial resolution.

\begin{table}[t]
\caption{\protect\label{ptentt} Number of iterations for
problem (\ref{ptent}), discretized with stepsizes (\ref{passi}) and
(\ref{Dt}), at each timestep $i\Delta t$, $i=1,\dots,20$.}
\centerline{\begin{tabular}{||r||r|r|r|r||}
\hline\hline
$N$ & 25 & 50 & 75 & 100\\ \hline \hline
$i\backslash n$ &  625 & 2500 & 5625 & 10000 \\ \hline \hline
$1$ &   5& 6& 8& 8\\
$2$ &   5& 6& 6& 6\\
\vdots &  \vdots& \vdots& \vdots& \vdots\\
$20$ &  5& 6& 6& 6\\
\hline\hline
\end{tabular}}
\end{table}

The second test problem is the evolutionary version of the elastic-plastic
torsion problem (\ref{epp}):
\begin{eqnarray}\nonumber
&&u_t \ge\Delta u + C, \quad u(x,y,t)\ge\psi(x,y),\quad (u_t-\Delta
u-C)(u-\psi)=0,\\
\label{pepp}\\ && (x,y,t)\in\Omega\times(0,\tau], \qquad
u\mid_{\partial\Omega} \equiv 0, \qquad u\mid_{t=0} =
\max\left(\psi,0\right)\equiv 0,\nonumber\end{eqnarray} where
$\psi$ and $\Omega$ are the same items defined in (\ref{epp}). The
spatial discretization is the same used for that problem (see
(\ref{passi1})), whereas the discretization in time is, for sake
of simplicity, done by means of the implicit Euler method, by
using a constant stepsize (\ref{Dt}). Also in this case, the
resulting PLS, to be solved at each time step, has dimension
$n=N^2$, whose matrix $T$ is the same as that obtained for problem
(\ref{epp}). Table~\ref{peppt} summarizes the obtained results, in
terms of required number of iterations, for $\tau=5$ and $\nu=20$.
In such a case, the approximation at $t=\tau$ is quite close to
the limit solution plotted in Figure~\ref{eppf}. The number of
iterations required for obtaining the solution turns out to be
quite similar to that listed in Table~\ref{eppt} for the
corresponding stationary problem.

\begin{table}[t]
\caption{\protect\label{peppt} Number of iterations for
problem (\ref{pepp}), discretized with stepsizes (\ref{passi1}) and
(\ref{Dt}), at each timestep $i\Delta t$, $i=1,\dots,20$.}
\centerline{\begin{tabular}{||r||rrrr|rrrr|rrrr|rrrr||}
\hline\hline
$N$ & \multicolumn{4}{c|}{25} & \multicolumn{4}{c|}{50} &
\multicolumn{4}{c|}{75} & \multicolumn{4}{c||}{100}\\ \hline \hline
$n$ &  \multicolumn{4}{c|}{625} & \multicolumn{4}{c|}{2500} &
\multicolumn{4}{c|}{5625} & \multicolumn{4}{c||}{10000} \\ \hline \hline
$i\backslash C$ &-5& -10&-15 &-20 & -5& -10&-15 &-20 & -5& -10&-15 &-20
& -5& -10&-15 &-20 \\
\hline\hline
$1$ &   9& 5& 4& 4& 17& 10& 7& 5& 25& 13& 9& 7& 32& 16& 11& 9\\
$\vdots$ &  $\vdots$& $\vdots$& $\vdots$& $\vdots$& $\vdots$& $\vdots$&
$\vdots$& $\vdots$& $\vdots$& $\vdots$& $\vdots$& $\vdots$& $\vdots$& $\vdots$&
$\vdots$& $\vdots$\\
$20$ &  9& 5& 4& 4& 17& 10& 7& 5& 25& 13& 9& 7& 32& 16& 11& 9\\
\hline\hline
\end{tabular}}
\end{table}

\section{Conclusions}\label{fine}

Two simple semi-iterative Newton-type procedures for solving
certain classes of piecewise linear systems have been
investigated. Such piecewise linear systems are derived from the
efficient modeling of obstacle problems. It has been shown that,
under rather general assumptions, the iterates are well defined
and monotonically converge to an exact solution of the given
system in a finite number of steps. A few numerical examples,
concerning both the classical obstacle problem, and its
evolutionary (i.e., parabolic) counterpart, prove the
effectiveness of the proposed methods.

\subsection*{Acknowledgements} The authors are indebted with Prof.
V.\,Casulli for his valuable comments.

\end{document}